\documentclass{aomart}

\fancypagestyle{firstpage}
{
  \fancyhf{}
  \cfoot{\footnotesize\thepage}
}

\usepackage{amsmath}
\usepackage{amssymb}
\usepackage[low-sup]{subdepth}

\newtheorem[{}\it]{lemma}{Lemma}[section]

\usepackage{color}
\usepackage{colortab}

\definecolor{blue}{rgb}{0.0,0.0,0.8}
\definecolor{gray1}{rgb}{0.80,0.80,0.80}
\definecolor{gray2}{rgb}{0.85,0.85,0.85}
\definecolor{gray3}{rgb}{0.90,0.90,0.90}
\definecolor{gray4}{rgb}{0.95,0.95,0.95}
\newcommand{\flag}{\textcolor{blue}}

\title[A short note on reduced residues]{A short note on reduced residues}
\author{Pascal Stumpf}

\begin{document}

\begin{abstract}
  We solve a problem due to Recam\'{a}n about the lower bound behavior 
  of the maximum possible length among all arithmetic progressions in the 
  least reduced residue system modulo $n$, as $n \to \infty$. 
\end{abstract}

\maketitle

\section{Introduction}

For any positive integer $n > 1$, let $\mathcal{A}(n) = \{a \in \mathbb{Z}^{+} : a < n, (a, n) = 1\}$ 
be the (nonempty) set of all smaller positive integers relatively prime to $n$, or 
in other words the least reduced residue system modulo $n$, and define $f(n)$ as 
the maximum possible length among all arithmetic progressions in $\mathcal{A}(n)$. In a 
letter from 1995 [1] Bernardo Recam\'{a}n asked if $f(n)$ tends to infinity with $n$, 
i.e. if for each $k \in \mathbb{Z}^{+}$ there exists a constant $n_{k}$ such that $\mathcal{A}(n)$ contains an 
arithmetic progression of length $k$ for all $n \geqslant n_{k}$. 

One very nice but deep result coming to mind here is that of Ben Green 
and Terence Tao [2] telling us about arbitrary long arithmetic progressions in 
the primes, and in fact it is a promising indicator for a positive answer to our 
question, since $\mathcal{A}(n)$ contains all primes less than $n$ except its prime factors. 
However, it turns out that we can prove the truth of our conjecture by using 
only elementary methods, and in what follows we like to present one possible 
(hopefully cute) solution. But before we start, let us consider a few examples 
to become even more familiar with our notations and the problem itself: 

\vspace{1.5ex}

\begin{center}
  \begin{tabular}{ccc}
    \LCC
    \color{gray1} & \color{gray2} & \color{gray1} \\
    $n$ & $\mathcal{A}(n)$ & $f(n)$ \\
    \ECC
    \LCC
    \color{gray3} & \color{gray4} & \color{gray3} \\
    $2$ & $\{\flag{1}\}$ & $1$ \\
    $3$ & $\{\flag{1}, \flag{2}\}$ & $2$ \\
    $4$ & $\{\flag{1}, \flag{3}\}$ & $2$ \\
    $5$ & $\{\flag{1}, \flag{2}, \flag{3}, \flag{4}\}$ & $4$ \\
    $6$ & $\{\flag{1}, \flag{5}\}$ & $2$ \\
    \ECC
    $\vdots$ & $\vdots$ & $\vdots$ \\
  \end{tabular}
  \quad
  \begin{tabular}{ccc}
    \LCC
    \color{gray1} & \color{gray2} & \color{gray1} \\
    $n$ & $\mathcal{A}(n)$ & $f(n)$ \\
    \ECC
    \LCC
    \color{gray3} & \color{gray4} & \color{gray3} \\
    $12$ & $\{\flag{1}, \flag{5}, 7, 11\}$ & $2$ \\
    $13$ & $\{\flag{1}, \flag{2}, \flag{3}, \flag{4}, \flag{5}, \flag{6}, \flag{7}, \flag{8}, \flag{9}, \flag{10}, \flag{11}, \flag{12}\}$ & $12$ \\
    $14$ & $\{\flag{1}, 3, \flag{5}, \flag{9}, 11, \flag{13}\}$ & $4$ \\
    $15$ & $\{\flag{1}, 2, \flag{4}, \flag{7}, 8, 11, 13, 14\}$ & $3$ \\
    $16$ & $\{\flag{1}, \flag{3}, \flag{5}, \flag{7}, \flag{9}, \flag{11}, \flag{13}, \flag{15}\}$ & $8$ \\
    \ECC
    $\vdots$ & $\vdots$ & $\vdots$ \\
  \end{tabular}
\end{center}

\newpage

\section{Ideas and Proof}

First, let us suppose $n$ is prime: then all of the numbers $1, \dots, n - 1$ are 
relatively prime to $n$ and form an arithmetic progression of length $n - 1$ with 
common difference $1$, which means $f(n) = n - 1$. In the more general case of 
a prime power $n = p^{r}$, where $p$ is prime and $r \in \mathbb{Z}^{+}$, we similarly still have 
$\{1, \dots, p - 1\} \subset \mathcal{A}(n)$ and thus $f(n) \geqslant p - 1$, but if $r \geqslant 2$ we can also look at 
the numbers $1 + m \cdot p$ for $0 \leqslant m < p^{r - 1}$, all of them lying in $\mathcal{A}(n)$ since none 
of them being divisible by $p$, and forming an arithmetic progression of length 
$p^{r - 1}$ with common difference $p$, giving us even $f(n) \geqslant p^{r - 1} = n / p$ here. 

Now let us consider squarefree numbers $n = p_{1} p_{2} \dots p_{d}$, where $d \geqslant 2$ and 
$2 \leqslant p_{1} < p_{2} < \ldots < p_{d}$ (odd) are prime. Like before, a good idea seems to be 
looking at numbers of the form $1 + m \cdot q$, this time choosing $q = p_{1} p_{2} \dots p_{d - 1}$ 
and $0 \leqslant m < p_{d}$, which ensures us that 
\[
  a_{m} = 1 + m \cdot q \leqslant 1 + (p_{d} - 1) \cdot q = 1 + n - q \leqslant 1 + n - 2 < n 
\]
is not divisible by any of the primes $p_{1}, p_{2}, \dots, p_{d - 1}$, although we are not sure 
about non-divisibility by $p_{d}$ yet. However, together $a_{0}, a_{1}, \dots, a_{p_{d} - 1}$ represent 
a complete residue system modulo $p_{d}$, because if $a_{x} \equiv a_{y} \hspace{-0.2ex}\pmod{p_{d}}$ for some 
$0 \leqslant x < y < p_{d}$ $(*)$, then $0 \equiv a_{y} - a_{x} = (y - x) \cdot q \hspace{-0.2ex}\pmod{p_{d}}$ and $(q, p_{d}) = 1$ 
would imply $(y - x) \equiv 0 \hspace{-0.2ex}\pmod{p_{d}} \Leftrightarrow x \equiv y \hspace{-0.2ex}\pmod{p_{d}}$ in contradiction to $(*)$. 
In particular only one member of $a_{0}, a_{1}, \dots, a_{p_{d} - 1}$ is divisible by $p_{d}$, say $a_{m}$, 
and so by the box principle we get that $a_{0}, \dots, a_{m - 1}$ or $a_{m + 1}, \dots, a_{p_{d} - 1}$ is an 
arithmetic progression of length at least $(p_{d} - 1) / 2$ with common difference $q$ 
completely contained in $\mathcal{A}(n)$, which delivers $f(n) \geqslant (p_{d} - 1) / 2$. 

Finally, let us introduce exponents $r_{1}, r_{2}, \dots, r_{d} \in \mathbb{Z}^{+}$ such that we can 
cover all remaining numbers $n = p_{1}^{r_{1}} p_{2}^{r_{2}} \dots p_{d}^{r_{d}}$, where $r_{1} + r_{2} + \ldots + r_{d} > d$. 
Because $n$ has the same prime factors as $p_{1} p_{2} \dots p_{d}$, we get $\mathcal{A}(p_{1} p_{2} \dots p_{d}) \subset 
\{a + m \cdot p_{1} p_{2} \dots p_{d} : a \in \mathcal{A}(p_{1} p_{2} \dots p_{d}), 0 \leqslant m < p_{1}^{r_{1} - 1} p_{2}^{r_{2} - 1} \dots p_{d}^{r_{d} - 1}\} = \mathcal{A}(n)$ 
by observing $(a, n) = 1 \Leftrightarrow (a, p_{1} p_{2} \dots p_{d}) = 1$ running over all integers $a$, and 
hence $f(n) \geqslant f(p_{1} p_{2} \dots p_{d}) \geqslant (p_{d} - 1) / 2$. On the other hand, we might again 
do better by looking at the numbers $1 + m \cdot p_{1} p_{2} \dots p_{d}$ forming an arithmetic 
progression of length $p_{1}^{r_{1} - 1} p_{2}^{r_{2} - 1} \dots p_{d}^{r_{d} - 1}$ with common difference $p_{1} p_{2} \dots p_{d}$, 
and both ideas in one lead us to $f(n) \geqslant \max \{(p_{d} - 1) / 2, n / (p_{1} p_{2} \dots p_{d})\}$. 

\vspace{1.5ex}

After we obtained lower bounds on $f(n)$ according to all possible prime 
factorizations of $n$, we are almost ready to prove our main result, but first let 
us collect them in the following more compact statement: 

\begin{lemma}
  For $n > 1$ we have $f(n) \geqslant \max \{(p - 1) / 2, n / P\}$, where $p$ is 
  the largest prime factor of $n$ and $P$ is the product of all prime factors of $n$. 
\end{lemma}

\newpage

\begin{lemma}
  For each $k \in \mathbb{Z}^{+}$ there exists a constant $n_{k}$ such that $\mathcal{A}(n)$ 
  contains an arithmetic progression of length $k$ for all $n \geqslant n_{k}$. 
\end{lemma}

\textit{Proof.} Let $P_{2k}$ be the product of all primes not exceeding $2k$ and define 
$n_{k} = k \cdot P_{2k} \geqslant 1 \cdot 2$. Moreover, let us fix any $n \geqslant n_{k}$ and (as in Lemma 2.1) 
denote its largest prime factor by $p$. If $p \geqslant 2k + 1$, we immediately arrive at 
$(p - 1) / 2 \geqslant ((2k + 1) - 1) / 2 = k$. But then in the other case $p < 2k + 1$, we 
note that all prime factors of $n$ do not exceed $2k$, implying their product $P$ 
divides $P_{2k}$, and so, in particular, $n / P \geqslant n_{k} / P = k \cdot P_{2k} / P \geqslant k$. Combining 
everything we get $f(n) \geqslant \max \{(p - 1) / 2, n / P\} \geqslant k$, and our claim follows. \hfill $\Box$ 

\vspace{1.5ex}

Captured by Lemma 2.2, we mainly worked on lower bounds so far and 
almost forgot about searching for possible upper bounds on $f(n)$. In order to 
catch up on them, let us change our point of view and conclude by showing: 

\begin{lemma}
  For $n > 1$ we have $f(n) \leqslant \max \{(p - 1) / 1, n / P\}$, where $p$ is 
  the largest prime factor of $n$ and $P$ is the product of all prime factors of $n$. 
\end{lemma}

\textit{Proof.} Suppose $a_{0}, a_{1}, \dots, a_{s - 1}$ is an arithmetic progression contained in 
$\mathcal{A}(n)$ with common difference $q$ and length $s$. Now we focus a bit more on $q$: 
If $q \geqslant P$, we can only come up to $s \leqslant n / P$, since otherwise $s > n / P$ implies 
$a_{s - 1} = a_{0} + (s - 1) \cdot q \geqslant 1 + ((n / P + 1) - 1) \cdot P \geqslant n + 1$ and our last member 
would not be in $\mathcal{A}(n)$ anymore. In the other case, we have $q < P$, yielding $q$ 
is missing at least one prime factor $p'$ of the squarefree number $P$ dividing $n$. 
But then $(q, p') = 1$ once again, like around $(*)$, whispers us that, whenever 
$s \geqslant p'$, the first $p'$ members $a_{0}, a_{1}, \dots, a_{p' - 1}$ do represent a complete residue 
system modulo $p'$, and thereby one of them, being a multiple of $p'$, could not 
ly within $\mathcal{A}(n)$ anymore, leaving us only $s \leqslant p' - 1 \leqslant p - 1$ left here. Uniting 
both cases we reach $f(n) \leqslant \max \{n / P, p - 1\}$, as desired. \hfill $\Box$ 

\vspace{1.5ex}

\textit{Acknowledgements.} All my thanks to Miriam (and Christian) for always 
encouraging me to write things down and making my first paper possible. 

\vspace{3.0ex}

\small

\begin{center}
  \textbf{References} 
\end{center}

\vspace{1.5ex}

\begin{center}
  \setlength{\tabcolsep}{2pt}
  \begin{tabular}{cp{352pt}}
    [1] & \textsc{Richard Guy}, \textit{Unsolved Problems in Number Theory}, Springer (2004), 146--147. 
    \\[0.0ex]
    [2] & \textsc{Ben Green} and \textsc{Terence Tao}, The primes contain arbitrarily long arithmetic 
          progressions, \textit{Annals of Mathematics} \textbf{167} (2008), 481--547. 
  \end{tabular}
\end{center}

\vspace{1.5ex}

\begin{center}
  \textit{e-mail:} \href{mailto:littlefriend@mathlino.org}{littlefriend@mathlino.org} 
\end{center}

\end{document}